\documentclass[a4paper,reqno]{amsart}
\usepackage{ 
	amsmath, 
	amssymb, 
	amsthm, 
	amscd, 
	amstext, 
	verbatim, 
	enumerate, 
	mathtools,
	color}

\usepackage{hyperref}
\usepackage{comment}
\usepackage[all]{xy}
\usepackage[utf8]{inputenc} %
\usepackage{lmodern}        %
\usepackage[T1]{fontenc}    %
\usepackage{hyperref}       %
\usepackage{url}            %
\usepackage{booktabs}       %
\usepackage{todonotes}
\usepackage{amsfonts}       %
\usepackage{nicefrac}       %
\usepackage{microtype}      %
\usepackage{xcolor}         %
\usepackage{mathtools}		%
\usepackage{aliascnt} 		%
\usepackage{enumitem}		%
\usepackage{graphicx}       %
\usepackage{caption}        %
\usepackage{subcaption}     %
\usepackage{float} 			%
\usepackage{forest}			%
\usepackage{mathrsfs} 		%
\setlist[enumerate,1]{label={(\roman*)}} %
\makeatletter
\newcommand{\myitem}[1]{%
	\item[#1]\protected@edef\@currentlabel{#1}%
}
\makeatother

\usepackage{tikz}
\usetikzlibrary{arrows.meta}
\usepackage{tikz-cd}

\hypersetup{
	colorlinks,
	linkcolor={blue!80!black},
	citecolor={blue!80!black},
	urlcolor={blue!80!black}
}

\theoremstyle{plain}
\newtheorem{theorem}{Theorem}[section]

\newaliascnt{propCt}{theorem}

\aliascntresetthe{propCt}

\newaliascnt{lemmaCt}{theorem}
\newtheorem{lemma}[lemmaCt]{Lemma}
\aliascntresetthe{lemmaCt}

\newaliascnt{corCt}{theorem}

\aliascntresetthe{corCt}

\theoremstyle{definition}
\newaliascnt{defiCt}{theorem}
\newtheorem{defi}[defiCt]{Definition}
\aliascntresetthe{defiCt}

\newaliascnt{convCt}{theorem}

\aliascntresetthe{convCt}

\newaliascnt{notaCt}{theorem}
\newtheorem{nota}[notaCt]{Notation}
\aliascntresetthe{notaCt}

\newaliascnt{remCt}{theorem}

\aliascntresetthe{remCt}

\newaliascnt{qstCt}{theorem}

\aliascntresetthe{qstCt}

\newaliascnt{pbmCt}{theorem}

\aliascntresetthe{pbmCt}

\newaliascnt{exaCt}{theorem}

\aliascntresetthe{exaCt}

\theoremstyle{plain}
\newtheorem{alphthm}{Theorem}			%

\newaliascnt{alphqstCt}{alphthm}

\aliascntresetthe{alphqstCt}

\newaliascnt{alphcorCt}{alphthm}

\aliascntresetthe{alphcorCt}

\newaliascnt{alphpropCt}{alphthm}

\aliascntresetthe{alphpropCt}

\newaliascnt{alphlemCt}{alphthm}
\newtheorem{alphlem}[alphlemCt]{Lemma}
\aliascntresetthe{alphlemCt}

\newaliascnt{alphproblemCt}{alphthm}

\aliascntresetthe{alphproblemCt}

\numberwithin{equation}{section}

\newcommand{\N}{\mathbb N}

\newcommand{\C}{\mathbb C}

\DeclareMathOperator{\Prob}{Prob}

\newcommand{\cint}[1]{\int_{\raisebox{-2pt}{$\scriptstyle #1$}}}

\begin{document}
\title[Topological and measurewise amenability]{Topological and measurewise amenability for locally compact groupoids}

\author[Julian Kranz]{Julian Kranz}
\address[Julian Kranz]{Universit\"at M\"unster, Mathematisches Institut, Einsteinstr. 62, 48149 M\"unster, Germany}
\email{julian.kranz@uni-muenster.de}
\urladdr{https://sites.google.com/view/juliankranz/}


\begin{abstract}
    We prove the equivalence of topological and measurewise amenability for $\sigma$-compact locally compact groupoids with a Hausdorff unit space and a continuous Haar system. %
	This in particular includes the second countable locally compact Hausdorff case.
	Our proof averages arbitrary Radon probability measures along strictly positive densities to obtain quasi-invariant measures. 
	We overcome the lack of absolute continuity of the original measure with respect to its quasi-invariant substitute by also averaging the Reiter functions along the same density.
	\begin{figure}[H]
		\includegraphics[width=.8\linewidth]{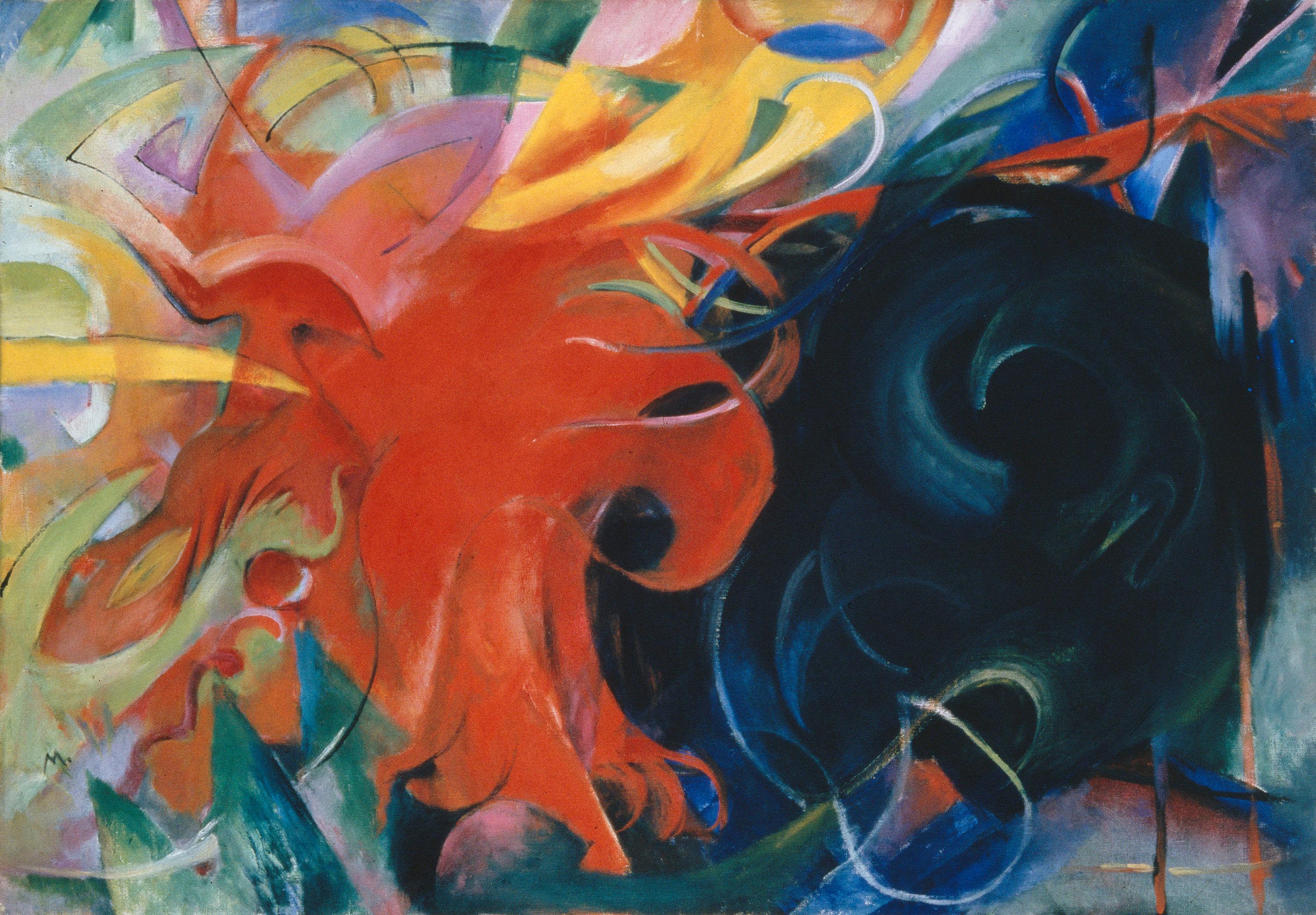}
		\caption*{{K\"ampfende Formen}, Franz Marc (1914)\protect\footnotemark}
	\end{figure}
\end{abstract} 
\maketitle
\footnotetext[1]{\url{https://commons.wikimedia.org/wiki/File:Franz_Marc_-_K\%C3\%A4mpfende_Formen_1914.jpg}, visited 10/28/2026.}
\tableofcontents

\section{Introduction}

Amenability was introduced as a group-theoretic averaging property in von Neumann's work on the Banach--Tarski paradox \cite{Neumann1929}.
It has since become a central approximation property in operator algebras and dynamics.
Its strength is illustrated by a plethora of classification results which turn abstract averaging properties into concrete internal approximations. 
These results can roughly be grouped into two highly intertwined categories: the \emph{measurable} and the \emph{topological} category.

In the measurable setting, the Connes--Feldman--Weiss theorem shows that countable amenable measured equivalence relations are hyperfinite \cite{Connes1981}.
The noncommutative analogues of measure spaces are von Neumann algebras, where a similar role is played by the Connes--Haagerup classification of injective factors with separable predual \cite{Connes1976,Haagerup1987}.

On the topological side, Connes' theorem has both inspired and technically facilitated the \emph{Elliott classification programme} which culminated in the classification of unital simple separable nuclear $\mathcal Z$-stable $C^*$-algebras satisfying the universal coefficient theorem (UCT) (see \cite{Winter2018,White2023} for an overview). 
Here, amenability not only enters through nuclearity, but it is also deeply connected to the UCT itself: By Tu's generalization of the Higson--Kasparov theorem \cite{Tu1999,Higson2001}, $C^*$-algebras associated to topologically amenable groupoids satisfy the UCT and their topological $K$-theory can be explicitly computed through methods from algebraic topology. 
The converse holds as well: Every unital simple separable nuclear $\mathcal Z$-stable $C^*$-algebra satisfying the UCT can be written as a twisted groupoid $C^*$-algebra \cite{Li2020}.

One reason for the power of amenability is its flexibility due to its many equivalent definitions.
While Reiter's definition \cite{Reiter1965} of amenable locally compact groups in terms of approximately invariant functions seemingly depends on the topology by referring to uniform convergence on compact sets, Zimmer's fixed point property \cite{Zimmer1978} only depends on the measurable structure of the group action. 
Both definitions of amenability were generalized to locally compact Hausdorff groupoids by Renault\footnote{See \cite{Adams1994} for the equivalence between Zimmer's definition and measurewise amenability.} \cite{Renault1980} and subsequently called \emph{topological} and \emph{measurewise} amenability in the foundational monograph \cite{AnantharamanDelaroche2000} by Anantharaman-Delaroche--Renault. Later, the theory was extended to non-Hausdorff groupoids in \cite{Renault2015}.
Since then it has been a long-standing open problem whether these two notions are equivalent (see \cite[page~92]{Renault1980}, \cite[\S~2.4]{AnantharamanDelaroche2023}).
Anantharaman-Delaroche and Renault proved the equivalence for locally compact Hausdorff 
groupoids with countable orbits, and more generally for groupoids equivalent to such groupoids \cite{AnantharamanDelaroche2000}.
This in particular includes all second countable Hausdorff \'etale groupoids. 
More recently, Buss--Echterhoff--Willett obtained the equivalence for transformation groupoids associated to actions of second countable locally compact groups on second countable locally compact Hausdorff spaces \cite{Buss2024}, building on \cite{Bearden2022}.
Our main result establishes the equivalence of topological and measurewise amenability for general $\sigma$-compact locally compact groupoids. 
\begin{alphthm}[\autoref{thm-main}]\label{thm-introduction}
Let $G$ be a $\sigma$-compact locally compact\footnote{Following Bourbaki's terminology, we use the word \emph{compact} to mean \emph{quasi-compact and Hausdorff}. In particular, locally compact spaces are assumed to be locally Hausdorff.} groupoid with a Hausdorff unit space and a Haar system.
Then $G$ is topologically amenable if and only if it is measurewise amenable.
\end{alphthm}
As a consequence, measurewise amenable $\sigma$-compact locally compact Hausdorff groupoids satisfy the Baum--Connes conjecture with coefficients and in the separable case their $C^*$-algebras satisfy the UCT \cite{Tu1999}. 

In order to explain the idea of the proof, we sketch the definitions of topological and measurewise amenability (see \autoref{sec-amenability} for the precise definitions) and recall Anantharaman-Delaroche--Renault's proof strategy. 
Following Reiter's condition \cite{Reiter1965}, topological amenability of a locally compact Hausdorff groupoid $G$ can be expressed in terms of the existence of compactly supported continuous functions $g_i\colon G\to [0,\infty)$ which are approximately invariant and approximately integrate to $1$ along range fibers \cite{Renault1980,AnantharamanDelaroche2000}. 
The approximations are quantified in terms of uniform convergence on compact subsets. 
Measurewise amenability is defined analogously, but the approximations are quantified in terms of weak-*-convergence with respect to \emph{quasi-invariant} Radon probability measures $\mu$ on the unit space $G^{(0)}$ of the groupoid. 
Here, quasi-invariance means that its pull-back $\mu\circ\lambda$ to $G$ along the Haar system is Radon--Nikodym-equivalent to its image under inversion. For transformation groupoids this simply means that $\mu$ is Radon--Nikodym-equivalent to all of its translates along the group action \cite{Renault1980}. 

It is easy to see that topological amenability implies measurewise amenability. The converse implication is usually proved using a Hahn--Banach argument: Since the \emph{Reiter functions} $(g_i)_{i}$ are closed under convex combinations, weak-*-approximations with respect to all Radon probability measures can be upgraded to uniform approximations on compact subsets. 
The remaining problem is: Not all Radon probability measures are quasi-invariant. 
If the groupoid has countable orbits, this problem can be overcome by averaging the original measure along the orbits. 
Since each point in the orbit is an atom, the original measure is absolutely continuous with respect to its quasi-invariant average. 
In particular, weak-*-approximations for the quasi-invariant measure yield weak-*-approximations for the original measure and the Hahn--Banach argument applies \cite{AnantharamanDelaroche2000}.

If the groupoid has uncountable orbits but is $\sigma$-compact, we can still average an arbitrary Radon probability measure $\mu$ along its orbits.
If $p\colon G\to (0,\infty)$ is a strictly positive probability density with respect to the Haar system, the convolved measure $p*\mu$ is quasi-invariant. 
In general, $\mu$ will fail to be absolutely continuous with respect to $p*\mu$. 
To overcome this problem, we simply average the Reiter functions as well.
Informally, this key observation can be phrased as the following averaging lemma:
\begin{alphlem}[\autoref{lem-averaging}]\label{lem-intro}
If $(g_n)_n$ is a Reiter sequence for $p*\mu$, then a subsequence of $(p*g_n)_n$ is a Reiter sequence for $\mu$.
\end{alphlem} 
Given the Reiter sequences from the above lemma, we follow a compactness argument from \cite{Renault2015} to construct a \emph{topological invariant mean} on $G$. Compared to the standard Hahn--Banach argument, this has the advantage of working in the non-Hausdorff case as well.

\subsection*{Acknowledgements}
This work was supported by the Deutsche Forschungsgemeinschaft (DFG, German Research Foundation) under Germany's Excellence Strategy EXC 2044/2-390685587, Mathematics M{\"u}nster: Dynamics--Geometry--Structure, and by the SFB 1442 of the DFG. 
The author would like to thank Siegfried Echterhoff for suggesting this topic during his PhD studies. The author would also like to thank Dennis Wulle for helpful comments and literature on Lie groupoids and Jamieson Bell for his help with AI tools.
\subsection*{Declaration of AI usage}
The author initially proved the main theorem for Hausdorff Lie groupoids, following the strategy of Ozawa--Suzuki \cite{Ozawa2021a}. 
Subsequently, he discussed the proof with a subscription-based generative AI model and discovered \autoref{lem-intro}. 
The author then reorganized the arguments, generalized the main theorem to $\sigma$-compact non-Hausdorff groupoids, and manually wrote the paper from scratch. 
Finally, he used generative AI for proofreading. 
While avoiding advertising AI companies publicly, the author is happy to disclose the used model upon request.
\section{Amenable groupoids}\label{sec-amenability}
We recall some basic definitions and fix notation.
By a \emph{compact} space we mean a Hausdorff topological space with the finite covering property.
A not necessarily Hausdorff space is called \emph{locally compact} if each point has a compact neighbourhood.
A topological space is called \emph{$\sigma$-compact} if it can be covered by countably many compact sets. 
For a locally compact Hausdorff space $X$, we denote by $C_c(X)$ the space of complex-valued compactly supported continuous functions on $X$. For a not necessarily Hausdorff locally compact space $X$, we denote by $\mathcal C_c(X)$ the linear span of all functions $f\colon X\to \C$ satisfying $f|_U\in C_c(U)$ and $f|_{X\setminus U}=0$ for some Hausdorff open subset $U\subset X$. 
We moreover denote by $\mathcal B(X)$ the space of all complex-valued Borel functions on $X$.
The subsets of non-negative functions are denoted by $\mathcal C_c(X)_+\subset \mathcal C_c(X)$ and $\mathcal B(X)_+\subset \mathcal B(X)$.

A \emph{groupoid} is a small category in which every morphism is invertible. Concretely, a groupoid consists of a set $G$ of morphisms, together with a subset $G^{(0)}\subset G$ of objects, range and source maps $r,s\colon G\to G^{(0)}$, a multiplication map $\cdot \colon G^{(2)}\to G$ defined on
\begin{equation}\label{eq-G2}
G^{(2)}\coloneqq \{(\gamma,\eta)\in G^2\mid s(\gamma)=r(\eta)\},
\end{equation}
and a (uniquely determined) inversion map $()^{-1}\colon G\to G$ satisfying a list of natural axioms which are summarized by the word ``category''. 
A \emph{topological groupoid} is a groupoid $G$ together with a topology such that the structure maps $r,s,\cdot$, and $()^{-1}$ are continuous. 
A \emph{locally compact groupoid} $G$ is a topological groupoid whose underlying topology is locally compact. 
We will frequently assume that the \emph{unit space} $G^{(0)}$ of $G$ is Hausdorff. In this case, the \emph{fibers} $G_x\coloneqq s^{-1}(x)$ and $G^x\coloneqq r^{-1}(x)$ are Hausdorff as well \cite[Proposition~2.8]{Tu2004}.
We refer to \cite{Paterson2012} for a textbook introduction to locally compact groupoids and to \cite{AnantharamanDelaroche2000,Williams2019} for an introduction to amenable groupoids.

As in the case of locally compact groups, the definition of amenability for locally compact groupoids requires the notion of an integral which is given by the \emph{Haar system}. The main difference to the group case is that for groupoids, its existence is not automatic (for instance, it forces the range map to be open \cite{Paterson2012}). 
\begin{defi}[Haar system]
Let $G$ be a locally compact groupoid with a Hausdorff unit space. A \emph{left Haar system} on $G$ is a family $\{\lambda^x\}_{x\in G^{(0)}}$ of Radon measures on $G^x$ of full support such that 
\begin{enumerate}
	\item For every $f\in \mathcal C_c(G)$, the map $\lambda(f)\colon G^{(0)}\to \C$ defined by
	\begin{equation}\label{eq-Haar}
	\lambda(f)(x)\coloneqq  \int f(\gamma)d\lambda^x(\gamma)
	\end{equation}
	is continuous,
	\item For every $\gamma\in G$ and every Borel subset $A\subset G^{s(\gamma)}$, we have $\lambda^{s(\gamma)}(A)=\lambda^{r(\gamma)}(\gamma A)$. 
\end{enumerate}
\end{defi}
If $\lambda$ is a Haar system, then any non-negative Borel function $f\colon G\to [0,\infty)$ defines a Borel function $\lambda(f)\colon G^{(0)}\to [0,\infty]$ by the same formula as in \eqref{eq-Haar} (see \cite{Renault2015}).
In order to abbreviate the definitions of topological and measurewise amenability, we introduce the following notation.

\begin{nota}
Let $G$ be a locally compact groupoid with a Hausdorff unit space and a Haar system $\lambda$.
For a non-negative Borel function $g\in \mathcal B(G)_+$, we write 
\[D_g\colon G\to [0,\infty],\quad D_g(\gamma)\coloneqq \int_{G^{r(\gamma)}}|g(\eta)-g(\gamma^{-1}\eta)|d\lambda^{r(\gamma)}(\eta).\]
\end{nota}
The following definitions are due to Renault \cite{Renault1980,AnantharamanDelaroche2000}.
\begin{defi}[Topological amenability]
	A $\sigma$-compact locally compact groupoid $G$ with a Hausdorff unit space and a Haar system $\lambda$ is \emph{topologically amenable} if there is a sequence of functions $(g_n)_{n\in \N}\subset \mathcal C_c(G)_+$ such that 
	\begin{enumerate}
		\item $\lambda(g_n) \leq 1$ for all $n\in \N$;
		\item $\lambda(g_n)\to 1$ uniformly on quasicompact subsets of $G^{(0)}$;
		\item $D_{g_n}\to 0$ uniformly on quasicompact subsets of $G$.
	\end{enumerate}
\end{defi}
As noted in \cite{Renault2015}, the use of sequences in the above definition is justified by our assumption of $\sigma$-compactness. 
\begin{defi}[Quasi-invariant measures]
	Let $G$ be a locally compact groupoid with a Hausdorff unit space and a Haar system $\lambda$. A \emph{quasi-invariant measure} on $G$ is a Radon probability measure $\mu$ on $G^{(0)}$ such that the measure $\mu\circ \lambda$ defined by 
	\[\int_Gf(\gamma)d(\mu\circ\lambda)(\gamma)=\int_{G^{(0)}}\int_{G^x}f(\gamma)d\lambda^x(\gamma)d\mu(x),\quad f\in \mathcal B(G)_+\]
	is Radon--Nikodym-equivalent to the measure $(\mu\circ\lambda)^{-1}$ defined by 
	\[\int_Gf(\gamma)d(\mu\circ\lambda)^{-1}(\gamma)=\int_Gf(\gamma^{-1})d(\mu\circ\lambda)(\gamma),\quad f\in \mathcal B(G)_+.\]
\end{defi}

The equivalence of the following definition to Renault's original definition in the second countable Hausdorff case is proved in \cite[Proposition~3.2.14]{AnantharamanDelaroche2000} (see also \cite[Theorem~10.49]{Williams2019}). 
In the $\sigma$-compact locally Hausdorff setting, we explicitly adopt it as our definition. 
\begin{defi}[Measurewise amenability]
	A $\sigma$-compact locally compact groupoid $G$ with a Hausdorff unit space and a Haar system $\lambda$ is called \emph{measurewise amenable} if for every quasi-invariant measure $\mu$ on $G^{(0)}$, there is a sequence $(g_n)_{n\in \N}$ of bounded Borel functions $g_n\colon G\to [0,\infty)$ such that 
	\begin{enumerate}
		\item $\lambda(g_n)\leq 1$ for all $n\in \N$;
		\item $\lambda(g_n)\to 1$ in the weak-*-topology of $L^\infty(G^{(0)},\mu)$;
		\item $D_{g_n}\to 0$ in the weak-*-topology of $L^\infty(G,\mu\circ \lambda)$.
	\end{enumerate}
\end{defi}

\section{Equivalence of the definitions}
The equivalence of topological and measurewise amenability for $\sigma$-compact locally compact groupoids relies on the fact that we can average both Radon probability measures on the unit space and Reiter functions on the groupoid along strictly positive Haar densities (see \autoref{lem-averaging}). 
The existence of a strictly positive Haar density is ensured by $\sigma$-compactness:
\begin{lemma}[Existence of positive Haar density]\label{lem-density}
	Let $G$ be a $\sigma$-compact locally compact groupoid with a Hausdorff unit space and a Haar system $\lambda$. Then there exists a Borel function $p\colon G\to (0,\infty)$ satisfying $\lambda(p)(x)=1$ for all $x\in G^{(0)}$. 
\end{lemma}
\begin{proof}
	Write $G=\bigcup_{n=1}^\infty K_n$ as a union of compact subsets $K_n\subset G$ and denote by $\mathbf 1_{K_n}$ their indicator functions. 
	Note that each $K_n$ is Borel by \cite[Proposition~3.6]{Niefield1983}. 
	Then 
	\[p_0\colon G\to (0,\infty),\quad p_0(\gamma)\coloneqq \sum_{n=1}^\infty\frac{\mathbf 1_{K_n}(\gamma)}{2^n(1+\lambda(\mathbf 1_{K_n})(r(\gamma)))}\] 
	defines a strictly positive function satisfying $0<\lambda(p_0)(x)<\infty$ for every $x\in G^{(0)}$. 
	Then
	\[p\colon G\to (0,\infty),\quad p(\gamma)\coloneqq \frac{p_0(\gamma)}{\lambda(p_0)(r(\gamma))}\]
	satisfies the desired properties. 
\end{proof}

The following lemma and its proof show that if $\mu$ is a not necessarily quasi-invariant Radon probability measure on $G^{(0)}$ and $(g_n)_{n\in \N}$ is a sequence of Reiter functions for the averaged measure $p*\mu=s_*(p(\mu\circ\lambda))$, then up to passing to subsequences, the averaged functions $(p*g_n)_{n\in \N}$ are a Reiter sequence for $\mu$. 
This is the key technical step which lets us adapt Anantharaman-Delaroche--Renault's argument for groupoids with countable orbits \cite{AnantharamanDelaroche2000} to $\sigma$-compact groupoids. 
\begin{lemma}[Averaging Reiter functions]\label{lem-averaging}
	Let $G$ be a measurewise amenable $\sigma$-compact locally compact groupoid with a Hausdorff unit space and a Haar system $\lambda$. 
	Let $\mu$ be a Radon probability measure on $G^{(0)}$. 
	Then there exists a sequence $(h_n)_{n\in \N}$ of Borel functions $h_n\colon G\to [0,\infty)$ and a $\mu$-conull subset $X_\mu\subset G^{(0)}$ satisfying 
	\begin{enumerate}
		\item\label{hn1} $\lambda(h_n)\leq 1$ for all $n\in \N$;
		\item\label{hn2} $\lambda(h_n)\to 1$ weak-$*$ in $L^\infty(G^{(0)},\mu)$;
		\item\label{hn3} $D_{h_n}(\gamma)\to 0$ for every $\gamma\in r^{-1}(X_\mu)$.
	\end{enumerate}
\end{lemma}
\begin{proof}
	Throughout the proof, we denote by $\mu\circ \lambda^2$ the measure on 
	\[G^{\times_r2}\coloneqq \{(\gamma,\eta)\in G\times G\mid r(\gamma)=r(\eta)\}\]
	defined by 
	\[\int_{G^{\times_r2}}\!\!\!f(\gamma,\eta)d(\mu\circ\lambda^2)(\gamma,\eta)\coloneqq \int_{G^{(0)}}\!\int_{G^x}\!\int_{G^x}\!\!\!f(\gamma,\eta)d\lambda^x(\gamma)d\lambda^x(\eta)d\mu(x),\quad f\in \mathcal B(G^{\times_r2})_+.\]
	Let $p$ be as in \autoref{lem-density}. 
	It follows from the proof of \cite[Proposition~2.9]{Renault2015} that the measure $\sigma\coloneqq p*\mu= s_*( p(\mu\circ \lambda))$ defined by 
	\[\int_{G^{(0)}}f(x)d\sigma(x)=\int_{G^{(0)}}\int_{G^x}f(s(\gamma))p(\gamma)d\lambda^x(\gamma)d\mu(x),\quad f\in \mathcal B(G^{(0)})_+\]
	is a Radon measure. We show that $\sigma$ is quasi-invariant. 
	For any $f\in \mathcal B(G)_+$, left-invariance of the Haar system implies
	\begin{equation}\label{eq-sigma-lambda}
		\begin{aligned}
			\int_Gf(\gamma)d(\sigma\circ \lambda)(\gamma)&=\int_{G^{(0)}}\left(\int_{G^x}f(\gamma)d\lambda^x(\gamma)\right)d\sigma(x)\\
		&=\int_{G^{(0)}}\int_{G^x}\left(\int_{G^{s(\eta)}}f(\gamma)d\lambda^{s(\eta)}(\gamma)\right)p(\eta)d\lambda^x(\eta)d\mu(x)\\
		&=\int_{G^{(0)}}\int_{G^x}\left(\int_{G^{x}}f(\eta^{-1}\gamma)d\lambda^{x}(\gamma)\right)p(\eta)d\lambda^x(\eta)d\mu(x)\\
		&=\int_{G^{\times_r2}}f(\eta^{-1}\gamma)p(\eta)d(\mu\circ \lambda^2)(\gamma,\eta).
		\end{aligned}
	\end{equation}
	Switching the roles of $\gamma$ and $\eta$ with $(\eta^{-1}\gamma)^{-1}=\gamma^{-1}\eta$ implies
	\begin{align*}
		\int_Gf(\gamma)d(\sigma\circ\lambda)^{-1}(\gamma)&=\int_Gf(\gamma^{-1})d(\sigma\circ \lambda)(\gamma)\\
		&=\int_{G^{\times_r2}}f(\eta^{-1}\gamma)p(\gamma)d(\mu\circ \lambda^2)(\gamma,\eta).
	\end{align*}
	Since both $(\gamma,\eta)\mapsto p(\gamma)$ and $(\gamma,\eta)\mapsto p(\eta)$ are strictly positive on $G^{\times_r2}$, we conclude that $\sigma\circ \lambda$ and $(\sigma\circ \lambda)^{-1}$ are Radon--Nikodym-equivalent measures. 
	Thus, $\sigma$ is quasi-invariant. 
	It moreover follows from $\lambda(p)=1$ that $\sigma$ is a probability measure.

	Since $G$ was assumed to be measurewise amenable, there is a sequence $(g_n)_{n\in \N}$ of bounded Borel functions $g_n\colon G\to [0,\infty)$ satisfying 
	\begin{enumerate}
		\item $\lambda(g_n)\leq 1$ for all $n\in \N$;
		\item $\lambda(g_n)\to 1$ weak-$*$ in $L^\infty(G^{(0)},\sigma)$;
		\item $D_{g_n}\to 0$ weak-$*$ in $L^\infty(G,\sigma\circ \lambda)$. 
	\end{enumerate}
	Using \eqref{eq-sigma-lambda}, the third property implies
	\[\int_Gp(\gamma)D_{g_n}(\gamma)d(\sigma\circ \lambda)(\gamma)=\int_{G^{\times_r2}}p(\eta^{-1}\gamma)D_{g_n}(\eta^{-1}\gamma)p(\eta)d(\mu\circ\lambda^2)(\eta,\gamma)\to 0.\]
	Since $p$ has positive values, we can thus pass to a subsequence of $(g_n)_{n\in \N}$ and find a $\mu$-conull subset $X_\mu\subset G^{(0)}$ satisfying 
	\begin{equation}\label{eq-subsequence}
		D_{g_n}(\eta^{-1}\gamma)\to 0\quad \text{ for } (\lambda^x\otimes \lambda^x) \text{ - a.\,e. }(\eta,\gamma)\in G^x\times G^x \text{ and } x\in X_\mu.
	\end{equation}
	We now claim that the functions
	\[h_n\colon G\to [0,\infty),\quad h_n(\gamma)\coloneqq (p*g_n)(\gamma)=\int_{G^{r(\gamma)}}p(\eta)g_n(\eta^{-1}\gamma)d\lambda^{r(\gamma)}(\eta)\]
	satisfy the desired conclusions.
	For any $x\in G^{(0)}$, left invariance of the Haar system implies 
	\begin{equation}\label{eq-lambda-hn}
	\begin{aligned}
		\lambda(h_n)(x)&=\int_{G^x}h_n(\gamma)d\lambda^x(\gamma)\\
		&=\int_{G^x}\int_{G^x}p(\eta)g_n(\eta^{-1}\gamma)d\lambda^x(\eta)d\lambda^{x}(\gamma)\\
		&=\int_{G^x}p(\eta)\int_{G^x}g_n(\eta^{-1}\gamma)d\lambda^{x}(\gamma)d\lambda^x(\eta)\\
		&=\int_{G^x}p(\eta)\int_{G^{s(\eta)}}g_n(\gamma)d\lambda^{s(\eta)}(\gamma)d\lambda^x(\eta)\\
		&=\int_{G^x}p(\eta)\lambda(g_n)({s(\eta)})d\lambda^x(\eta)\\
		&\leq \int_{G^x}p(\eta)d\lambda^x(\eta)\\
		&\leq 1,
	\end{aligned}
\end{equation}
verifying property \ref{hn1}.
For property \ref{hn2}, note that \eqref{eq-lambda-hn} implies
\[\|1-\lambda(h_n)\|_{L^1(G^{(0)},\mu)}=1-\int_{G^{(0)}}\lambda(h_n)(x)d\mu(x)=1-\int_{G^{(0)}}\lambda(g_n)(x)d\sigma(x)\to 0.\]
Since $\mu$ is a probability measure and $(1-\lambda(h_n))_{n\in \N}$ is bounded, this implies $\lambda(h_n)\to 1$ in the weak-$*$-topology of $L^\infty(G^{(0)},\mu)$. 
We proceed with verifying property \ref{hn3}.
Fix $\gamma\in r^{-1}(X_\mu)$.
Using left invariance of the Haar system at the third and sixth step, the equality $\int_{G^{r(\gamma)}} (p(\xi)-p(\gamma^{-1}\xi))d\lambda^{r(\gamma)}(\xi)=0$ at the fourth equality, and $\int_{G^{r(\gamma)}}p(\zeta)d\lambda^{r(\gamma)}(\zeta)=1$ at the fifth equality, we obtain 
\begin{align*}
	D&_{h_n}(\gamma)=\cint{G^{r(\gamma)}}\!\!\!\!\!\!|h_n(\eta)-h_n(\gamma^{\scriptscriptstyle-1}\eta)|d\lambda^{\scriptscriptstyle{r(\gamma)}}(\eta)\\
	=&\cint{{G^{\scriptscriptstyle{r(\gamma)}}}}\!\left|\cint{G^{\scriptscriptstyle{r(\gamma)}}}\!\!\!\!p(\xi)g_n(\xi^{\scriptscriptstyle{-1}}\eta)d\lambda^{\scriptscriptstyle{r(\gamma)}}\!(\xi) - \cint{G^{s(\gamma)}}\!\!\!\!p(\xi)g_n(\xi^{\scriptscriptstyle{-1}}\gamma^{\scriptscriptstyle{-1}}\eta)d\lambda^{s(\gamma)}\!(\xi)   \right|d\lambda^{\scriptscriptstyle{r(\gamma)}}\!(\eta)\\
	=&\cint{G^{\scriptscriptstyle{r(\gamma)}}}\!\left|\cint{G^{\scriptscriptstyle{r(\gamma)}}}\!\!\!\!(p(\xi)-p(\gamma^{\scriptscriptstyle{-1}}\xi))g_n(\xi^{\scriptscriptstyle{-1}}\eta)d\lambda^{\scriptscriptstyle{r(\gamma)}}\!(\xi) \right|d\lambda^{\scriptscriptstyle{r(\gamma)}}\!(\eta)\\
	=&\cint{G^{\scriptscriptstyle{r(\gamma)}}}\!\left|\cint{G^{\scriptscriptstyle{r(\gamma)}}}\!\!\!\!\!\!\!(p(\xi)-p(\gamma^{\scriptscriptstyle{-1}}\xi))\!\left(\!g_n(\xi^{\scriptscriptstyle{-1}}\eta)-\!\cint{G^{\scriptscriptstyle{r(\gamma)}}}\!\!\!\!\!\!\!\!\!p(\zeta)g_n(\zeta^{\scriptscriptstyle{-1}}\eta)d\lambda^{\scriptscriptstyle{r(\gamma)}}\!(\zeta)\!\right)\!d\lambda^{\scriptscriptstyle{r(\gamma)}}\!(\xi) \right|d\lambda^{\scriptscriptstyle{r(\gamma)}}\!(\eta)\\
	=& \cint{G^{\scriptscriptstyle{r(\gamma)}}}\!\left|\cint{G^{\scriptscriptstyle{r(\gamma)}}}\!\cint{G^{\scriptscriptstyle{r(\gamma)}}}\!\!\!\!\!\!\!\!\!\!\!p(\zeta)(p(\xi)-p(\gamma^{\scriptscriptstyle{-1}}\xi))\!\left(g_n(\xi^{\scriptscriptstyle{-1}}\eta)\!-\!g_n(\zeta^{\scriptscriptstyle{-1}}\eta)\right)\!d\lambda^{\scriptscriptstyle{r(\gamma)}}\!(\xi)d\lambda^{\scriptscriptstyle{r(\gamma)}}\!(\zeta)\right|\! d\lambda^{\scriptscriptstyle{r(\gamma)}}\!(\eta)\\
	\leq&  \cint{G^{\scriptscriptstyle{r(\gamma)}}}\!\cint{G^{\scriptscriptstyle{r(\gamma)}}}\!\!\!\!\!\!\!|p(\xi)-p(\gamma^{\scriptscriptstyle{-1}}\xi)|p(\zeta)\left(\cint{G^{\scriptscriptstyle{s(\zeta)}}}\!\!\!\!\!\!\!|g_n(\eta)-g_n(\xi^{\scriptscriptstyle{-1}}\zeta\eta)|d\lambda^{s(\zeta)}\!(\eta)\right)d\lambda^{\scriptscriptstyle{r(\gamma)}}\!(\zeta) d\lambda^{\scriptscriptstyle{r(\gamma)}}\!(\xi)\\
	=&\cint{G^{\scriptscriptstyle{r(\gamma)}}}\!\cint{G^{\scriptscriptstyle{r(\gamma)}}}\!|p(\xi)-p(\gamma^{\scriptscriptstyle{-1}}\xi)|p(\zeta)D_{g_n}(\zeta^{\scriptscriptstyle{-1}}\xi)d\lambda^{\scriptscriptstyle{r(\gamma)}}\!(\zeta) d\lambda^{\scriptscriptstyle{r(\gamma)}}\!(\xi).
\end{align*}
The right-hand integrand is bounded above by $2(p(\xi)+p(\gamma^{\scriptscriptstyle{-1}}\xi))p(\zeta)$ since 
\begin{equation}\label{eq-dominated}
	D_{g_n}(\zeta^{-1}\xi)\leq \lambda(g_n)(s(\zeta))+\lambda(g_n)(s(\xi))\leq 2.
\end{equation}
It thus follows from the dominated convergence theorem and \eqref{eq-subsequence} that $D_{h_n}(\gamma)\to 0$ for all $\gamma\in r^{-1}(X_\mu)$. 
\end{proof}

From \autoref{lem-averaging}, we now follow standard techniques for proving the equivalence of topological and measurewise amenability. 
In the Hausdorff case, we could conclude topological amenability with a convex Hahn--Banach argument. 
Since the non-Hausdorff case causes certain technical complications with the space of compactly supported continuous functions, 
we instead use a variant of the definition of topological amenability in terms of topological invariant means \cite{AnantharamanDelaroche2000}. 
Let $G$ be a $\sigma$-compact locally compact groupoid with a Hausdorff unit space and a Haar system $\lambda$. We denote by $\mathcal E=C_0(G^{(0)},L^1(G,\lambda))$ the completion of $\mathcal C_c(G)$ with respect to the norm 
\[\|f\|_\mathcal E\coloneqq \sup_{x\in G^{(0)}}\int_{G^x}|f(\gamma)|d\lambda^x(\gamma),\quad f\in \mathcal C_c(G).\]
The Banach space dual $\mathcal E^*$ can be identified with all complex Radon measures on $G$ that are of the form $\nu=\varphi(\mu\circ \lambda)$ for some Radon probability measure $\mu$ on $G^{(0)}$ and $\varphi\in L^\infty(G,\mu\circ\lambda)$ (see \cite[Proposition~2.9]{Renault2015}). 
The subspace of all $\varphi(\mu\circ \lambda)\in \mathcal E^*$ for a fixed $\mu$ is denoted by $M(G,\mu\circ\lambda)$ so that 
\begin{equation}\label{eq-MGL}
	\mathcal E^*=\bigcup_{\mu\in \Prob(G^{(0)})}M(G,\mu\circ\lambda). 
\end{equation}
We write $\nu\geq 0$ if $\nu$ is a positive measure, i.e. if $\nu(f)\geq 0$ for all $f\in \mathcal C_c(G)_+$. 
For a function $f\in \mathcal C_c(G)$, the convolution operator
\[f*-\colon\mathcal C_c(G)\to \mathcal C_c(G),\quad f*g(\gamma) \coloneqq\int_{G^{r(\gamma)}}f(\eta)g(\eta^{-1}\gamma)d\lambda^{r(\gamma)}(\eta),\quad g\in \mathcal C_c(G),\,\gamma\in G\]
extends continuously to $\mathcal E$ and thus to $\mathcal E^{**}$ by bitransposition. 
\begin{defi}[Topological invariant mean]\label{defi-topo-mean}
	Let $G$ be a $\sigma$-compact locally compact groupoid with a Hausdorff unit space and a Haar system $\lambda$. A \emph{topological invariant mean} is an element $m\in \mathcal E^{**}$ satisfying the following properties. 
	\begin{enumerate}
		\item $m(\nu)\geq 0$ for all $\nu\in \mathcal E^*$ with $\nu\geq 0$;
		\item $\|m\|\leq 1$;
		\item $m(\mu\circ\lambda)=1$ for all Radon probability measures $\mu$ on $G^{(0)}$;
		\item $f*m=(\lambda(f)\circ r)\cdot m$ for all $f\in \mathcal C_c(G)$. 
	\end{enumerate}
\end{defi}
\begin{theorem}[\cite{AnantharamanDelaroche2000,Renault2015}]\label{thm-top-mean}
A $\sigma$-compact locally compact groupoid with a Hausdorff unit space and a Haar system is topologically amenable if and only if it admits a topological invariant mean. 	
\end{theorem}

\begin{theorem}\label{thm-main}
	Let $G$ be a $\sigma$-compact locally compact groupoid with a Hausdorff unit space and a Haar system. Then $G$ is measurewise amenable if and only if $G$ is topologically amenable.
\end{theorem}
\begin{proof}
	We only prove the forward implication since the converse is classical \cite{Renault1980}.
	Our proof closely follows the proof of the implication $(i)\Rightarrow(ii)$ of \cite[Theorem~2.14]{Renault2015}.
	Assume that $G$ is measurewise amenable.
	For a finite set $P\subset \Prob(G^{(0)})$ of Radon probability measures on $G^{(0)}$  we denote by $\mathcal M_{P}\subset \mathcal E^{**}$ the set of elements $m\in \mathcal E^{**}$ satisfying 
	\begin{enumerate}[label=(M\arabic*)]
		\item\label{m1} $m(\nu)\geq 0$ for all $\nu\in \mathcal E^*$ with $\nu\geq 0$;
		\item\label{m2} $\|m\|\leq 1$;
		\item\label{m3} $m(\rho\circ\lambda)=1$ for all $\rho \in P$;
		\item\label{m4} $(f*m)|_{M(G,\rho\circ\lambda)}=((\lambda(f)\circ r)\cdot m)|_{M(G,\rho\circ\lambda)}$ for all $f\in \mathcal C_c(G)$ and $\rho\in P$,
	\end{enumerate}
	where $M(G,\rho\circ\lambda)$ is defined as in the paragraph above \eqref{eq-MGL}.
	Note that each $\mathcal M_P$ is bounded and weak-*-closed and thus weak-*-compact by the Banach--Alaoglu theorem.
	Moreover, we have $\mathcal M_P\subset \mathcal M_Q$ whenever $Q\subset P$. 
	By the finite intersection property, it thus suffices to show that each $\mathcal M_{P}$ is non-empty. 
	Indeed, in this case, any element of the intersection \[\bigcap_{P\underset{\text{finite}}\subset \Prob(G^{(0)})}\mathcal M_{P}\not=\emptyset\] will be a topological invariant mean in the sense of \autoref{defi-topo-mean} and $G$ will be topologically amenable by \autoref{thm-top-mean}. 
	Fix a finite subset $\emptyset \not=P\subset \Prob(G^{(0)})$. 
	We apply \autoref{lem-averaging} to $\mu\coloneqq \frac 1 {|P|}\sum_{\rho \in P}\rho\in \Prob(G^{(0)})$ and obtain a sequence $(h_n)_{n\in \N}$ of Borel functions $h_n\colon G\to [0,\infty)$ and a $\mu$-conull set $X_\mu\subset G^{(0)}$ satisfying 
	\begin{enumerate}[label=(H\arabic*)]
		\item\label{hn1b} $\lambda(h_n)\leq 1$ for all $n\in \N$;
		\item\label{hn2b} $\lambda(h_n)\to 1$ weak-$*$ in $L^\infty(G^{(0)},\mu)$;
		\item\label{hn3b} $D_{h_n}(\gamma)\to 0$ for every $\gamma\in r^{-1}(X_\mu)$.
	\end{enumerate}
	As in \cite[Theorem~2.14]{Renault2015}, we define elements $m_{h_n}\in \mathcal E^{**},\,n\in \N$ by 
	\[m_{h_n}(\varphi (\rho\circ \lambda))\coloneqq \int_Gh_n(\gamma)\varphi(\gamma)d(\rho\circ\lambda)(\gamma),\quad \varphi(\rho\circ\lambda)\in \mathcal E^{*}.\]
	Pass to a weak-*-convergent subnet $(m_{h_{n_i}})_{i\in I}$ of $(m_{h_n})_{n\in \N}$ and write $m\coloneqq \lim_{i\in I}m_{h_{n_i}}$. 
	We prove that $m$ satisfies properties \ref{m1} - \ref{m4} above. 
	Property \ref{m1} is immediate from the fact that each $h_n$ takes nonnegative values. 
	Property \ref{m2} follows from property \ref{hn1b}. 
	Property \ref{m3} follows from \ref{hn2b} and the fact that every $\rho\in P$ is absolutely continuous with respect to $\mu$. 
	In order to prove property \ref{m4}, fix $f\in \mathcal C_c(G)$ and $\varphi(\rho\circ \lambda)\in M(G,\rho\circ\lambda)$ for some $\rho\in P$. 
	Since $\rho\leq |P|\mu$, we can find $\psi\in L^\infty(G,\mu\circ \lambda)$ such that $\psi(\mu\circ \lambda)=\varphi(\rho\circ\lambda)$. 
	Using Fubini's theorem at the third step and the triangle inequality at the last step, we obtain
	\begin{align*}
		&\left|\left((f*m_{h_n})-(\lambda(f)\circ r)\cdot m_{h_n} \right)(\varphi(\rho\circ\lambda))\right|\\
		&=\left|\int_G((f*h_n)(\gamma)-\lambda(f)(r(\gamma))h_n(\gamma))\psi(\gamma)d(\mu\circ \lambda)(\gamma)\right|\\
		&=\left|\int_{G^{(0)}} \int_{G^x} \int_{G^{r(\gamma)}}\left(f(\eta)h_n(\eta^{-1}\gamma)-f(\eta)h_n(\gamma)\right)d\lambda^{r(\gamma)}(\eta) \psi(\gamma)d\lambda^x(\gamma)d\mu(x)\right|\\
		&=\left|\int_{G^{(0)}} \int_{G^x} \int_{G^{x}}f(\eta)(h_n(\eta^{-1}\gamma)-h_n(\gamma)) \psi(\gamma)d\lambda^x(\gamma)d\lambda^{x}(\eta)d\mu(x)\right|\\
		&\leq\|\psi\|_{\infty}\int_{G^{(0)}} \int_{G^x}|f(\eta)|D_{h_n}(\eta)d\lambda^{x}(\eta) d\mu(x).
	\end{align*}
By the same argument as in \eqref{eq-dominated}, the integrand is bounded by $2|f|$. 
Thus, \ref{hn3b} and the dominated convergence theorem imply that 
\[\left|\left((f*m_{h_n})-(\lambda(f)\circ r)\cdot m_{h_n} \right)(\varphi(\rho\circ\lambda))\right|\xrightarrow{n\to \infty}0.\] 
Since $m$ is the weak-*-limit of $(m_{h_{n_i}})_{i\in I}$, this proves property \ref{m4}.
\end{proof}
\bibliography{references.bib}
\bibliographystyle{alphaurl}

\end{document}